\documentclass[12pt]{article}
\usepackage{amssymb,latexsym, amsmath, amscd, array, graphicx}
\usepackage{amssymb}

\usepackage{setspace}
\usepackage{enumitem}
\input xy
\xyoption{all}
\makeatletter 
\usepackage{color}
\usepackage{comment}

\usepackage{setspace}
\usepackage{enumitem}

\date{}

\newtheorem{theorem}{Theorem}[section]

\newtheorem{proposition}[theorem]{Proposition}

\newtheorem{remark}[theorem]{Remark}
\newtheorem{problem}[theorem]{Problem}

\newcommand{\z}{{\Bbb Z}}

\newcommand{\re}{{\Bbb R}}
\newcommand{\N}{{\Bbb N}}
\newcommand{\s}{{\Bbb S}}
\newcommand{\C}{{\Bbb C}}

\newcommand{\lo}{\longrightarrow}

\newcommand{\black}{{\blacksquare}}

\newcommand{\mdim}{{\rm mdim}}
\newcommand{\mesh}{{\rm mesh}}
\newcommand{\ord}{{\rm ord}}
\begin{document}

\title {\bf A free $\z$-action by isometries on a compact metric space which is not embeddable into a cubical shift}
\author{ Alexander Dranishnikov 
\footnote{A.D.  is thankful to IMPAN (Warsaw), FIM ETH (Zurich), MPI fur Mathematik (Bonn) for hospitality
during his visits in 2024 when this project was initiated.} and Michael Levin\footnote{
 M.L. 
 was partially supported by the Israel Science Foundation grants No. 2196/20 and  3241/24,  and
 IMPAN (Warsaw)  during his visit in 2025.}
}
\maketitle
\begin{abstract} We construct an example announced in the title. It answers in a strong way a well-known open problem 
in topological dynamics.  {In fact our construction is an existence theorem}.  It is based on 
a Borsuk-Ulam type theorem whose proof heavily relies on $p$-adic completions of $G$-complexes and
the equivariant Sullivan conjecture.
\\
{\bf Math. Subj. Class.:}  37B02  (55P60, 55P91)

\end{abstract}

\begin{section}{Introduction}
\label{introduction}
The cubical shift $([0,1]^k)^\z $ is a dynamical system  on    $([0,1]^k)^\z$  with the action of $\z$ generated  by
the left shift 
$$(\sigma(x))_i =x_{i+1},\ \ 
x=(x_i) \in ([0,1]^k)^\z,\  i \in \z.$$

In 1974 Jaworski  \cite{jaworski} proved   that a free dynamical system $(X, \z)$ on 
a finite dimensional compact metric space $X$  is (equivariantly) embeddable  into  the shift $[0,1]^\z$. 
The natural question whether the finite dimensional restriction  on  $X$ is essential in Jaworski's theorem
  stayed unsettled for quite a while and 
only in 2000  was affirmatively answered by  Lindenstrauss  and Weiss  \cite{lindenstrauss-weiss}.
 Lindenstrauss-Weiss' result   is based on the notion of  mean dimension,
a topological invariant of dynamical systems introduced by 
 Gromov in 1999.
 Gromov's mean-dimension $\mdim  X $ of a dynamical system  $(X,\z)$ 
is  defined and  discussed in Section \ref{mean-dimension}, here we just mention that $\mdim X$ is
 a non-negative real number or $\infty$ and 
among basic properties of mean dimension are that $\mdim  ([0,1]^k)^\z=k$ and  $\mdim X  \leq k$ if $X$ 
is embeddable into $([0,1]^k)^\z$.  
 Thus  a minimal  $\z$-action on a compact metric {space} $X$
with $\mdim X > 1$ constructed by Lindenstrauss and Weiss \cite{lindenstrauss-weiss} is not embeddable
into $[0,1]^\z$ and {it} shows that Jaworski's theorem does not extend to infinite dimensional spaces.

The next break-through result  came  also in 2000. {It} is  due to  Lindenstrauss \cite{lindenstrauss} and  says that a minimal $\z$-action
on a compact metric space  $X$  with finite mean dimension, $\mdim X <\infty$, is embeddable into a cubical shift $([0,1]^k)^\z$
for some $k$.  {Also} Lindenstrauss
provided an estimate for $k$ that depends  on $\mdim X$  and  which was later
improved by Gutman and Tsukamoto \cite{gutman-tsukamoto} to an optimal one, see Section \ref{mean-dimension}
for more details.

The following natural generalization of Lindenstrauss' embedding theorem is considered as one of the major open
problems concerning the embeddability  into cubical shifts.
\begin{problem}
\label{problem}
Is {every} free $\z$-action on a compact metric space  with finite mean dimension  embeddable into a cubical  shift?
\end{problem}

It was widely believed  that the problem  has a positive answer  and for quite a while it was not even known whether there {are}
dynamical systems  satisfying the assumptions of Problem \ref{problem} and  not covered by the existing embedding theorems. 
However in 2022 Tsukamoto, Tsutaya and Yoshinaga  in 
their nice paper \cite{tsukamoto+} followed by Shi's clever refinement  \cite{shi}  constructed such  examples and, as the result,  showed
that a positive answer to Problem \ref{problem} is out of reach by  the  currently available embedding  techniques. 
But  still there was a hope that   known approaches may   provide a positive  answer to  Problem \ref{problem} in the important case of
mean dimension $0$. A more detailed description of the related results
is given in Section \ref{mean-dimension}. The goal of this paper is
to give a negative answer to Problem \ref{problem} even in mean dimension $0$.

\

{\bf Theorem {A}.}
{\em There is a free $\z$-action by isometries on a compact metric space which is not embeddable in
a cubical shift.}

\

A dynamical system with a $\z$-action by isometries  not only has  mean dimension $0$
 but also topological entropy $0$ (see Section \ref{mean-dimension}). 
 Thus Theorem A negatively answers Problem \ref{problem} in quite a strong way.
 Theorem A is derived from the following Borsuk-Ulam type theorem.
 
 Let  $\s^{2n-1}$  be  the unit sphere in $\C^{n}$. By the standard
 action of the $m$-cyclic group $\z/m\z$ on $\s^{2n-1}$ we mean the action
  $$gz =e^{2g\pi i/m}(z_1, \dots, z_n)$$ where $ g \in \z/m\z$ and $z=(z_1,\dots,z_n)\in \s^{2n-1}\subset \C^n$.
  Clearly $\z/m\z$ acts on $\s^n$ freely and  by isometries.
  
   Given a prime number $p$ and a sequence of  natural numbers $n_i, i\geq 1,$
  for each $i$ we consider
  the  sphere $\s^{2n_i-1}$ with the standard  free action of the group
 $\z/p^i \z$ and let  $$T_i : \s^{2n_i-1}\lo \s^{2n_i-1}$$ be the action of a generator of  $\z/p^i \z$. {We} denote by $X$ 
 the infinite product of all the spheres $\s^{2n_i-1}, i\geq 1,$ with the  diagonal $\z$-action generated by
 the homeomorphism  
 $ T : X \lo X$ which is defined by $$(T(x))_i=T_i (x_i)\ \ \text{for}\ \ x=(x_i)\in X
  \ \text{with}\ \ x_i \in \s^{2n_i-1}.$$
 Clearly  $\z$ acts on $X$ freely and by isometries. We will refer to $(X, \z)$ as the dynamical system 
 defined by the prime $p$ and the sequence $n_i$. 
 
  For a map $f : X \lo Y$ of the above $\z$-space $X$ to a topological space $Y$ we define 
 the  map to the product $f^\z : X \lo Y^\z$ by the formula $$(f^\z(x))_j=f(T^j(x)),\ \ {j \in \z},\ \ x \in X.$$
 We refer to $f^\z$  as
 the map induced by $f$ and note  that $f^\z$ is a $\z$-map with respect to the shift action of $\z$ on $Y^\z$.
 Note that $f$ coincides with the map $f^\z$  followed  by the projection of $Y^\z$ to  the $0$-coordinate factor $Y$.

\
 
 {\bf Theorem {B}.}
 {\em For every  prime number $p$  there is a sequence  of   natural numbers  $n_i$  such that
 the dynamical system $(X,\z)$ defined by $p$ and the sequence $n_i$ has the property:
  For every natural $k$
 and every map $f : X \lo \re^k$ there is a point $x \in X$ such that  $f (\z x)\subset \re^k$ contains at most $p^{k-1}$ points.
 
 Moreover,
 the set  $f^\z(\z x) \subset (\re^k)^\z$  
 contains at most $p^{k-1}$  points.
 
 Here  $\z x$ is the orbit of $x$ under the $\z$-action on $X$
 and  $f^\z : X \lo (\re^k)^\z$ is the map 
 induced  by $f$.}
 
 \

Borsuk-Ulam type theorems  have  a long history and  a wide variety of methods used  for their proofs. However,
despite a considerable effort, 
  the authors of this paper were not able to  prove Theorem B
  applying known approaches.
 The approach that finally led to  proving  Theorem B seems to be new and 
 heavily relies on $p$-adic completions of $G$-complexes and the equivariant Sullivan conjecture.
 \\
 
  Theorem A trivially follows from Theorem B.
 \\
 {\bf Proof of  Theorem A}. We will show that  a dynamical system $(X,\z)$ satisfying the conclusions 
 of Theorem B is not embeddable  in a cubical shift.  Aiming at a contradiction assume
 that there is an equivariant embedding $\phi : X \lo ([0,1]^k)^\z$  and let $f: X \lo [0,1]^k$ be the map $\phi$ 
 followed by the projection of $([0,1]^k)^\z$ to the $0$-coordinate cube  $[0,1]^k$. 
 Note that $\phi=f^\z$  and, hence,
  $f^\z(\z x)$ is infinite  for every $x \in X$. This
 contradicts   the conclusion of Theorem B that there is $x$ in $X$ such that
 $f^\z(\z x)$ has at most $p^{k-1}$ points. $\black$
 \\
 
 In Section \ref{section-borsuk-ulam}  we present  Theorem \ref{generalization-2} which generalizes
 Theorem B to compact group actions  where the groups are  countable products  of 
 finite $p$-groups. This broader framework allows us to show (in  Theorem \ref{generalization-1})
  that  the group  $\z$  in Theorem A
 can be replaced by any infinite subgroup of a countable product of finite $p$-groups.
 \\
 
 The paper is organized as follows: In Section \ref{mean-dimension} we provide
 a general background related to mean dimension and embedding theorems, 
 in Section \ref{reduction} we reduce Theorem B to  Proposition \ref{proposition-for-borsuk-ulam}
 which involves  $p$-adic completions of $G$-spaces,
 Sections  \ref{section-sullivan} and  \ref{g-maps} are  devoted 
 to auxiliary results needed for Proposition \ref{proposition-for-borsuk-ulam},
 in Section \ref{section-borsuk-ulam} we prove  Proposition \ref{proposition-for-borsuk-ulam}
  and present generalizations of  Theorems  A and  B.

\end{section}

\begin{section}{Mean dimension and embedding theorems}\label{mean-dimension}
Let $(X, \z)$ be a dynamical system on a compact metric $X$ generated by a homeomorphism
$T : X \lo X$. For a cover $\cal U$ of $X$ we denote by $\mesh(\cal U)$ the supremum of the diameters
of the sets in $\cal U$
and by  $\ord(\cal U)$ the maximal number of sets in $\cal U$ whose intersection is non-empty.

 The mean dimension  $\mdim X$   of $(X,\z)$   can be defined in the following short way.
Let $d$ be  a positive real number. We will use the notation $\mdim X \leq d^-$ to say that
for every $\epsilon >0$ there is an open cover $\cal U$ of $X$ such that
$\mesh (T^i ({\cal U}))< \epsilon$ for every $i$ satisfying $0\leq id \leq \ord (\cal U)$.
Then Gromov's mean dimension coincides with 
$$\mdim X=\inf \{d : \mdim  X \leq d^-\}.$$

In a similar way one can define the topological entropy $h_{top}  X$ of $(X, \z)$. 
We will write $h_{top} X \leq d^-$ if for every $\epsilon>0$ there is an open cover
$\cal U$ of $X$ such that $\mesh({T^i}\cal U) < \epsilon$ for every $i$
satisfying $0\leq id \leq \log |{\cal U}| $ where $|{\cal U}|$ is the number of elements of $\cal U$.
Then $$h_{top} X =\inf \{d : h_{top} X \leq d^- \}.$$

Since $\mesh (T^i {\cal U})=\mesh ({\cal U})$ for every $i$  if $T$ is an isometry, we immediately get
that $\mdim X=h_{top} X=0$ if  $\z$ acts by isometries on $X$. It is known
that for every $\z$ action 
$\mdim X=0$ if $h_{top} X $ is finite, $h_{top} X < \infty$  \cite{lindenstrauss-weiss}.
 Lindenstrauss \cite{lindenstrauss}
showed that the latter relation is even deeper:  {E}very minimal dynamic system $(X,\z)$ with
$\mdim X =0$ is the inverse limit of dynamical systems with finite topological entropy.

A crucial role in
Lindenstrauss' embedding theorem \cite{lindenstrauss}  is played by so-called level functions.
Gutman \cite{marker-finite-dimension} singled out a property that allows to construct level functions and coined
this property  as the marker property. 

A dynamical system $(X, \z)$ generated by 
a homeomorphism $T: X \lo X$ has the {\em marker property }
if for every non-negative integer $n$ there is an open subset $U$ of $X$ such that 
$\z U =X$ and the sets $U, TU,\dots, T^n U$ are disjoint. 

 Lindenstrauss' embedding theorem  actually says 
that there is a constant $C >0$ (namely $C=1/36$)  such that 
a dynamical system $(X,\z)$ having the marker property embeds in the shift $([0,1]^k)^\z$ if  $\mdim X < Ck$. 
In a  striking work~\cite{gutman-tsukamoto}  Gutman and Tsukamoto  showed using 
a {signal processing} technique  that the constant
$C$   can be improved to $C=1/2$ and in an earlier work Lindenstrauss and Tsukamoto 
\cite{lindenstrauss-tsukamoto} showed that the constant $C=1/2$ is best possible.
Let us emphasize that all known  constructions of embeddings of infinite dimensional dynamical systems 
into cubical shifts  essentially rely on the marker property.

It was  conjectured for  a while that  any free dynamical system has the marker property. A result supporting 
this conjecture was a theorem by Gutman \cite{marker-finite-dimension} saying that a free $\z$-action 
on a finite dimensional compact metric space  has the marker property. Tsukamoto, Tsutaya and Yoshinaga   \cite{tsukamoto+} 
disproved 
this conjecture in $\mdim=\infty$  and later Shi \cite{shi} refined \cite{tsukamoto+} by constructing free dynamical systems $(X, \z)$
of arbitrarily small $\mdim$ and not having the marker property. 

However it was still conjectured \cite{tsukamoto+} that
a free dynamical system $(X,\z)$ with $\mdim=0$ has the marker property. 
Theorem  A clearly disproves this conjecture as well.
Moreover,  the dynamical system $(X, \z)$ from 
 Theorem B satisfies  a condition which is stronger than just not having 
 the marker property.  Namely, for any open subsets $U_0$ and $U_1$
 of $X$ such that 
 $\z U_0=\z U_1=X$ we have $U_0 \cap U_1 \neq \emptyset$.
Indeed, assume that  $U_0 \cap U_1 = \emptyset$. Then,
 since $X $ is compact,  there are closed subsets $F_0 \subset U_0$ and 
$F_1 \subset U_1$ such that $\z F_0 =\z F_1 =X$. 
Take a map $f : X\lo [0,1]$ sending $F_0$ to $0$ and $F_1$ to $1$. Then for every $ x \in X$
the set $f(\z x)$ contains at least two points, {namely} the end points $\{ 0, 1\}$, which contradicts Theorem B for $k=1$.

\end{section}

\begin{section}{{Theorem B and equivariant $p$-adic completions}}
\label{reduction}
Let $G$ be a discrete group and $X$ a $G$-space i.e., a topological space with an  action  of $G$.
For a subgroup $H$ of $G$ we denote $X^H=\{ x \in X: hx=x, h \in H\}$, the fixed point set of the action of $H$ on $X$.
By a $G$-map between  $G$-spaces we mean a $G$-equivariant map.
We consider $X \times [0,1]$ as a  $G$-space with 
 the product action of the action  of $G$ on $X$ and the trivial action on $[0,1]$.  This   defines in the obvious way 
 the notions of  $G$-homotopy of $G$-maps and $G$-homotopy equivalence of  $G$-spaces.
For a topological space $Y$ the (full) $G$-shift is a  $G$-space $Y^G$ with the product topology and 
the shift action of $ G$ defined for $h \in G$ 
by  $$(hy)_g =y_{gh}\ \ \text{where}\ \ y=(y_g) \in Y^G\ \ \text{and}\ \ g \in G.$$
For a map $f : X \lo Y$ the {\em induced $G$-map} $f^G : X \lo Y^G$ is defined by
$$(f^G(x))_g=f(gx), \ \ x \in X\ \ \text{and}\ \ g\in G.$$

We say that a $G$-space $K$ is a $G$-CW-complex if $K$ is a  CW-complex,
 $G$ acts on $K$ by sending cells to cells and if a cell  $C$ of $K$ is invariant under
the action of  an element  $ g \in G$ then  $g$ acts trivially on $C$.
By a simplicial $G$-complex we mean a $G$-CW-complex whose  cell structure comes
from a triangulation. A free simplicial $G$-complex means a simplicial $G$-complex on which $G$ acts freely.
Two  $G$-spaces are isomorphic if there is a $G$-map between them which is a homeomorphism.

Let $p$  be a prime number and  $G_m=\z/p^m\z$, $m\geq 1$. For a natural number $l \geq 1$ 
we denote  $E(l, m)=(\re^l)^{G_m}$ and consider
$E(l,m)$ with the shift  action of $G_m$ on $(\re^l)^{G_m}$. For $0\leq q < m$ we denote 
by $E^q(l,m)$ the set of $p^q$-periodic points of $E(l,m)$, i.e., the set of all points of $E(l,m)$ whose orbit size
is at most $p^q$,  and   set $$E(l,m,q)=E(l,m)\setminus E^q(l,m).$$ 
Clearly, $E(l,m,q)$ is invariant 
under the action of $G_m$ on $E(l,m)$. Thus, we consider $E(l,m,q)$ with the action of $G_m$. 

{We regard} $E(l,m)$ as a Euclidean space
 and $E^q(l,m)$ as its linear subspace. Let 
$K(l,m,q)$ be the unit sphere in the orthogonal complement  $E^q(l,m)^\perp $ of $E^q(l,m)$  in $E(l,m)$. Note 
that $K(l,m,q)$ is invariant under the action of $G_m$ on $E(l,m)$ and we consider $K(l,m,q)$ with  the action
of $G_m$. {For $0\leq q<m$  we denote by $K^m(l,m+1,q)$ the set of $p^m$-periodic points of $K(l,m+1,q)$.}

 \begin{proposition}\label{general-properties}
 ${}$
 
 \begin{enumerate}[label=\upshape(\Roman*),ref=\theproposition (\Roman*)]
     \item\label{general-1}  $K(l,m,q)$ and $E(l, m,q)$ are $G_m$-homotopy equivalent.
  \item\label{general-2}
  Let $l \leq l'$. Then $K(l,m,q)$ admits a $G_m$-map to  $K(l',m,q)$.
  \item\label{general-3}
    $K^m(l, m+1, q)$ 
     with the action of the quotient group $G_m=G_{m+1}/H$, $H=\z/p\z,$ is isomorphic to $K(l,m,q)$.
\end{enumerate}     

\end{proposition}
{\bf Proof.} \\
{\bf (I). } The orthogonal projection
of $E(l,m)$ to $E^q(l,m)^\perp$  restricted to $E(l,m,q)$ and followed by the radial projection 
$$E^q(l,m)^\perp \setminus \{ 0 \} \lo K(l,m,q)$$ is a $G_m$-homotopy equivalence between $E(l,m,q)$
and  $K(l,m,q)$.
  \\\\
{\bf (II).}  Consider an embedding $e : \re^{l} \lo \re^{l'}$ of 
$\re^{l}$ into $\re^{l'}$ as a coordinate subspace. Then $e$ induces
the natural linear  isometric $G_m$-embedding $e_m : E(l,m)\lo E(l', m)$ such that
$$e_m(E^q(l,m)^\perp)=E^q(l',m)^\perp \cap e_m(E(l,m)).$$ This implies that
$e_m$  also $G_m$-embeds  $K(l,m,q)$ into $K(l',m,q)$.
\\
\\ 
{\bf (III).} Consider the diagonal linear embedding $$e : E(l,m)=({\re^l)}^{p^m} \lo E(l, m+1)=(({\re^l})^{p^m})^{p}$$
defined by $$e(x)=C(x, \dots ,x),\ x \in ({\re^l})^{p^m},$$   where $(x,\dots, x)$ is a $p$-tuple and
 the real constant $C >0$ is chosen such that $e$ is an isometric embedding.
Note that  $e$  sends $E(l,m)$ onto  $E^m(l,m+1)$ with $e(E^q(l,m))=E^q(l,m+1)$  {for all $q$} and $e$ is  an isomorphism onto  $E^m(l,m+1)$ considered
with the action of $G_{m+1}/H$. This implies that $e$  isomorphically sends $K(l,m,q)$ onto $K^m(l,m+1,q)$ considered with the action
of $G_{m+1}/H$. $\black$
\\

The next proposition presents  an important  property of $K(l,m,q)$. 
\begin{proposition}
\label{crucial}
 Let $X$ be a  $G_m$-space such that 
      $X$ does not admit a $G_m$-map to $K(l,m,q),$ $ 0\leq q< m$. Then  for every 
     map $f : X \lo \re^{l}$ there is a point $x\in X$ such that $f^{G_m}(G_m x)$ contains 
  at most    $p^{q}$ points  where  $f^{G_m} : X \lo (\re^l)^{G_m}$  is
  the $G_m$-map induced by $f$. 
  
  In particular, $f(G_m x)$ contains at most $p^q$ points.
  \end{proposition}
  {\bf Proof.} By Proposition \ref{general-1},  $K(l,m,q)$  and $E(l,m,q)=E(l,{m}) \setminus E^q(l,m)$ are
  $G_m$-homotopy equivalent. We recall  that  $$E(l,m)=(\re^l)^{G_m}\ \ \text{and}\ \ E(l,m,q)=E(l,m)\setminus E^q(l,m).$$
  Thus, {we have a $G_m$-map} $f^{G_m}: X\lo E(l,m)$ {where} $X$ does not admit a $G_m$-map to  $E(l,m,q)$.
  Then there is $x \in X$ such that $f^{G_m}(x) \in E^q(l,m)$. This means that 
   $f^{G_m}(x)$ is $p^{q}$-periodic and hence
  $f^{G_m}(G_m x)=G_m f^{G_m}(x)$ contains at most $p^{q}$ points. $\black$
  \\

 Given a sequence of natural numbers $n_i$, $i\geq 1,$ consider  the sphere $\s^{2n_i-1}$
 endowed with the standard action of $G_i=\z/p^i\z$ and let 
   $$N_m{=\prod_{i=1}^{i=m}\s^{2n_i-1} }$$ be  the product of spheres $\s^{2n_i-1}$, $1\leq i \leq m$, endowed  with the diagonal action of
  $G_m$. Clearly $N_m$ admits a triangulation for which  $N_m$ is
  a free simplicilal $G_m$-complex. We refer to $N_m$ as the $G_m$-complex determined by $p$ and the sequence $\{n_i\}$.
   Theorem B is derived from
 
 \begin{proposition}
 \label{proposition-for-borsuk-ulam-simplified}
 For every prime $p$ and every  sequence of natural numbers $l_q$, $q\geq 0,$ 
 there is a sequence of  natural numbers $n_i$, $i\geq 1,$
 such for  each $m\geq 1$ the  $G_m$-complex $N_m$ determined by $p$ and  the sequence $n_i$ 
   does not admit a $G_m$-map to
  $K(l_q,m,q)$ for every $0\leq q < m$. 
\end{proposition}

Proposition \ref{proposition-for-borsuk-ulam-simplified} implies even a  more general version of 
Theorem B.

 \begin{theorem} \label{general-borsuk-ulam}
 For every prime $p$ and every sequence $l_q$, $q\geq 0,$ of natural numbers 
  there is a sequence  of   natural numbers  $n_i, i\geq 1,$  such that
 the dynamical system $(X,\z)$ defined by $p$ and the sequence $n_i$ (see Section \ref{introduction})
  has the property:
  For every  $q\geq 0$
 and every map $f : X \lo \re^{l_q}$ there is a point $x \in X$ such that  the set  $f^\z(\z x)$  
 contains at most $p^{q}$  points where the $\z$-map $f^\z : X \lo (\re^{l_q})^\z$ is induced by $f$.
 
 In particular, $f(\z x)$ contains at most $p^{q}$ points.
 \end{theorem}
{\bf Proof.} We will show that   the sequence 
$n_i$ satisfying  the conclusions of Proposition \ref{proposition-for-borsuk-ulam-simplified} 
for the sequence $l_q$ will have the required properties.  

Suppose this is not true and let    a map $f : X \lo  \re^{l_q}$ be  such that 
 $f^\z(\z x)$ contains more than $p^{q}$ points  for every $x \in X$.
 Then, since $f^\z$ is continuous and $X$ is compact, there is $\epsilon>0$ such that
 for every $x \in X$ the set $f^\z(\z x)$ contains $p^q+1$ points 
 with pairwise distances larger than $\epsilon$. 
 
 Recall that $X$ is the product of  all the spheres $\s^{2n_i-1}, i\geq 1,$ and
 $N_m$ is the product of the spheres $\s^{2n_i -1}$ for $1 \leq i \leq m$.
  Let $\phi_m : X \lo N_m$ be the projection.
 Then for   a sufficiently large $m$ 
 there is  a map $f_m : N_m \lo \re^{l_q}$ such that $(f_m \circ \phi_m)^\z$ and $f^\z$ are
 $\epsilon/2$-close, and   hence $(f_m \circ \phi_m)^\z (\z x)$ contains at least 
 $p^q+1$ points for every $x \in X$. Clearly we may assume that $m >q$.
 {We} regard the action of $G_m$ on $N_m$  as  the action  of $\z$ whose isotropy group is $p^m \z$
 and note that  with respect to this action  $\phi_m$ is a $\z$-map.
 Thus we  formally consider $N_m$ with the actions of $G_m$ and $\z$ that gives rise to 
 the induced maps $$f_m^{G_m}: N_m \lo (\re^{l_q})^{G_m}\ \ \text{and}\ \  f_m^{\z}: N_m \lo (\re^{l_q})^\z$$
 for which we have $$|f_m^{G_m}(G_m y)|=|f_m^\z (\z y)|$$ for all $y \in N_m$.
 Also for all $x \in X$ we have  $$( f_m \circ \phi_m)^\z (\z x)= f_m^\z (\phi_m (\z x))=
 f^\z_m(G_m\phi_m(x)).$$ Thus we get that for all $x \in X$  $$|f_m^{G_m}(G_m \phi_m(x))|\ge p^q+1$$
 and, since $\phi_m$ is surjective  and $m>q$, we
{obtain}  a contradiction with
Proposition \ref{crucial}.  $\black$
\\
\\
{\bf Proof of Theorem B}: follows from  Theorem \ref{general-borsuk-ulam} for the sequence $l_q=q+1, q\geq 0$. 
$\black$
\\

The way we prove Proposition \ref{proposition-for-borsuk-ulam-simplified} is by proving 
a stronger version, namely replacing $K(l_q, m, q)$ by its $p$-adic completion.

The $p$-adic completion  $K^\land_p$ of a CW-complex  $K$  is discussed in detail in Section \ref{section-sullivan}.
What is  important for us  now is  that the $p$-adic completion is a functor from the category of  CW-complexes 
and continuous maps  to itself. The $p$-adic completion of a map $f : K\lo {L}$ of CW-complexes
is denoted by
$$f^\land_p : K^\land_p \lo {L}^\land_p.$$
If $K$ is a $G$-space then $K^\land_p$ automatically turns into
a $G$-space with $g \in G$  acting on $K^\land_p$ by the $p$-adic completion of the action of $g$ on $K$.
Then for a $G$-map $f : K \lo {L}$ we get that $f^\land_p : K^\land_p \lo {L}^\land_p$ is a $G$-map as well.
An important property of $p$-adic completion  is that the $p$-adic completion of a $G$-CW-complex is again a $G$-CW-complex. 
Thus we get that the $p$-adic completion is also a functor from the category of $G$-CW-complexes and $G$-maps to itself.
Another important property of $p$-adic completion  is
that  for a  $G$-CW-complex $K$ there is a {naturtal} $G$-map from $K$ to $K^\land_p$ (Proposition \ref{basic-0}). 

Thus  we reduce  Proposition \ref{proposition-for-borsuk-ulam-simplified}  to

 \begin{proposition}
 \label{proposition-for-borsuk-ulam}
 For every prime $p$ and every sequence of natural numbers $l_q$, $q\geq 0,$ there is a sequence of  natural numbers $n_i$, $i \geq 1,$
 such for  each $m\geq 1$ the $G_m$-complex $N_m$ determined by $p$ and  the sequence $n_i$ 
   does not admit a $G_m$-map to
  $K(l_q,m,q)^\land_p$ for {each}  $0\leq q < m$. 
\end{proposition}
Proposition \ref{proposition-for-borsuk-ulam-simplified} immediately follows from  Proposition \ref{proposition-for-borsuk-ulam}.
\\
{\bf Proof of Proposition \ref{proposition-for-borsuk-ulam-simplified}.}
 As we already  mentioned 
 $K(l_q, m, q)$ admits a $G_m$-map to $K(l_q, m, q)^\land_p$. Hence, if  
  $N_m$  does not admit a $G_m$-map to $K(l_q, m, q)^\land_p$  then
 $N_m$ does not admit a $G_m$-map to $K(l_q, m, q)$  either. $\black$
 \\
 \\
Proposition \ref{proposition-for-borsuk-ulam} is proved in Section \ref{section-borsuk-ulam}, the  proof is short and the reader
can easily follow  it right now relying on very few {results explicitly stated in Sections \ref{section-sullivan} and \ref{g-maps}}.
\begin{remark}
\label{remark}
${}$
\end{remark} 
 {We} note that the only case when $K(l,m,q)$ is not connected is  $p=2$, $l=1, m=1, q=0$ and
then it is a  $2$ point set $\s^0$. The only case when $K(l,m,q)$ is connected but not simply connected
is $p=3,l=1,m=1,q=0$ and then it is a circle $\s^1$.
 By  Proposition \ref{general-2},   $K(l,m,q)$ admits a $G_m$-map to $ K(l',m,q)$  for $l\leq l'$. Then, 
 since the $p$-adic completion is functorial, we get that 
 $K(l,m,q)^\land_p$ admits a $G_m$-map to $ K(l',m,q)^\land_p$.
 This clearly implies that  we can replace in Proposition \ref{proposition-for-borsuk-ulam} 
(and also in Proposition \ref{proposition-for-borsuk-ulam-simplified}) the sequence $l_q$ by any
sequence $l'_q$ with $l_q \leq l'_q$  and  assume without loss of generality 
that  $K(l_q, m, q)$ is connected  (and even simply connected
if we wish) for every prime $p$.

\end{section}

\begin{section}{$p$-Adic Completion}
\label{section-sullivan} All the groups are assumed to be discrete and $p$  is a prime number.
 For a simplicial set $X$, $|X|$ stands for  the CW-complex which is
 the   geometric realization of $X$.
Let $K$ be a CW-complex. By ${\rm Sin} K$ we denote the singular simplicial set of $K$ and 
by the $p$-adic completion of $K$  we mean the CW-complex  $$K^\land_p=|R_\infty  {\rm Sin} K|$$
where $R_\infty$ is the Bousfield-Kan completion functor with $R=\z/p\z$ \cite{b-k}.  Note that  for every simplicial set $X$
there is a natural simplicial map $X \lo R_\infty X$.
Also note that the adjunction map $|{\rm Sin} K| \lo K$ is a homotopy equivalence. Therefore there is a canonical up to homotopy map $j_K:K\to K^\land_p$.

Now let $G$  be a  discrete group
and $K$ be a $G$-CW-complex. 
 This turns in a functorial way ${\rm Sin} K$ into a simplicial $G$-set and 
$|{\rm Sin} K|$ and $K^\land_p$ into $G$-CW-complexes. Again  the functoriality of the realization and singular functors
implies that for
every subgroup $H$ of $G$ 
the map $|{\rm Sin} K| \lo K$  restricted to 
$$|{\rm Sin} K|^H=|({\rm Sin} K)^H|=|{\rm Sin} (K^H)| \lo K^H$$  
is a homotopy equivalence and hence, by the equivariant Whitehead theorem,
$|{\rm Sin} K| \lo K$ is a $G$-homotopy equivalence.

For $G$-CW-complexes $K$ and $L$ we denote by ${\rm Map}_G (L,K)$ the space of $G$-maps 
from $L$ to $K$ endowed with the compact-open topology.
By  $EG$ we denote  a contractible  $G$-CW-complex with a free action of $G$.
{Thus, $EG$ is the  universal cover of an Eilenberg-MacLane complex $K(G,1)$.}

\begin{proposition}
\label{basic}
Let $G$ be a  group and  $K$ a $G$-CW-complex.
${}$

 \begin{enumerate}[label=\upshape(\Roman*),ref=\theproposition (\Roman*)]
     \item\label{basic-0} A canonical map $j_K:K\to K^\land_p$ can be taken to be a  $G$-map.
    \item\label{basic-1} 
     Let $G$ be a $p$-group, 
    $H$  a normal subgroup of $G$, and $K$ a finite $G$-CW-complex.
    Then   $(K^H)^\land_p$ and $ (K^\land_p)^H$ are $G$-homotopy equivalent.
\end{enumerate}
\end{proposition}
{\bf Proof.}
 Property \ref{basic-0} follows from the fact that the adjunction map $|{\rm Sin} K| \lo K$ is a $G$-homotopy equivalence.

Property \ref{basic-1} can be derived  from \cite{carlsson} as follows. Let $A$ be  a subgroup of $G$.
 Note that $K^A$ and $(K^\land_p)^A$ are 
the fixed point sets of  $K$ and $K^\land_p$ respectively considered as  $A$-CW-complexes.
Also note  $$({\rm Sin} K)^A\ \ \text{and}\ \ (R_\infty {\rm Sin} K)^A,\ \ R=\z/p\z,$$ are the fixed point sets
of ${\rm Sin} K$ and $R_\infty {\rm Sin} K$ respectively considered as simplicial $A$-sets.
Denote by $X={\rm Sin} K$. Then, by Proposition II.8  of \cite{carlsson},  the inclusion $X^A\lo  X$
induces a weak homotopy equivalence  $$R_\infty (X^A)\lo (R_\infty X)^A.$$ 
This implies that the corresponding  map 
$$|R_\infty (X^A)| \lo |(R_\infty X)^A|$$ is  a weak homotopy equivalence.
Note that ${\rm Sin} (K^A)=({\rm Sin} K)^A=X^A$ and hence 
 $$(K^A)^\land_p =|R_\infty{\rm Sin}(K^A)|=|R_\infty (X^A)|.$$ Also note
that  $$|(R_\infty X)^A|=|R_\infty X|^A=(K^\land_p)^A$$
 and we get  that the map $$(K^A)^\land_p \lo (K^\land_p)^A$$ is a weak homotopy equivalence.
 
 Thus, by replacing $A$ by $AH$,  we get that the map 
$$(K^{AH})^\land_p \lo (K^\land_p)^{AH}$$  
  is a weak homotopy equivalence  and, replacing $K$ by $K^H$, we also get that 
  the map   $$((K^H)^A)^\land_p  \lo ((K^H)^\land_p)^A$$  is a weak homotopy equivalence.
  Note that $(K^H)^A=K^{AH}$ and  $((K^\land_p)^H)^A=(K^\land_p)^{AH}$
  and  hence  the map $$((K^H)^\land_p)^A \lo ((K^\land_p)^H)^A$$ is a weak homotopy equivalence.
  Since, it holds for every subgroup $A$ of $G$  and the map   $$(K^H)^\land_p \lo (K^\land_p)^H$$ is a $G$-map,
  the equariant Whitehead theorem implies that 
  $(K^H)^\land_p$ and $(K^\land_p)^H$ are $G$-homotopy equivalent. $\black$
\\

\begin{theorem}[Equivariant Sullivan Conjecture]~\label{EqSuCo}
 Let $G$ be a $p$-group and  $K$ a finite $G$-CW-complex. Then,
    after identifying the constant maps in ${\rm Map}_G   (EG, K^\land_p)$ with $(K^\land_p)^G$,
    the inclusion  $$(K^\land_p)^G\to {\rm Map}_G(EG, K^\land_p)$$  is a weak homotopy equivalence.
\end{theorem}
The equivariant Sullivan conjecture {also known as the generalized Sullivan conjecture} was proved independently 
by Carlsson~\cite{carlsson},  Dwyer–Miller–Neisendorfer~\cite{d-m-n}, and  Lannes \cite{lannes}.

The following proposition  is derived from \cite{paper-b-k, b-k}.
\begin{proposition}
\label{tower}
Let $K$ be a  connected $G$-CW-complex whose homology groups with coefficients in $\z/p\z$  are finite. Then
there is a tower of fibrations $\{\Omega_s\}$ with  bonding maps $\omega_{s} : \Omega_{s} \lo \Omega_{s-1}$
such  that  $\Omega_s$  is a connected $G$-CW-complex, $\Omega_s$ has finite homotopy groups, 
$\omega_s$ is a fibration which is also a $G$-map,  and there is a $G$-map $K^\land_p
\to \Omega=\varprojlim \Omega_s$ which is a weak homotopy equivalence.
\end{proposition}
{\bf Proof.} Let  $X={\rm Sin} K$ and let $\{ R_s X\}$ be the Bousfield-Kan tower of fibrations for $X$ and $R=\z/p\z$
\cite{b-k}. 
We recall that  $R_\infty X=\varprojlim R_sX$.
Since the homology groups $H_*(K; R)$ are finite,
the homology groups $H_*(X; R)$ are finite as well and hence all groups $E_r^{s,t}(X; R)$ of 
 the Bousfield-Kan homotopy spectral sequence of $X$  with coefficients in $R$
 are finite,   Ch. VI, 9.3 \cite{b-k}. 
The extended homotopy spectral sequence of the tower $\{ R_s  X\}$ coincides
 in dimensions$\geq 1$  with the  spectral sequence  $E_r^{s,t}(X;R)$,  Ch. I, 4.4 \cite{b-k} and
 \cite{paper-b-k}. Then examining the $E_1$-term of  the Bousfield-Kan homotopy spectral sequence  
 we conclude that the homotopy groups of $R_s X$ are finite. Thus, the extended homotopy spectral sequence is Mittag-Leffler 
and converges to the homotopy groups of the inverse limit $\pi_*(R_\infty X)=\pi_*(|R_\infty X|)$.

Recall that the  geometric realization of a (Kan) fibration is a (Serre) fibration.
 Denote $\Omega_s =| R_s  X|$  and  let the map $\omega_{s} : \Omega_{s} \lo \Omega_{s-1}$ be
 induced by the fibration  $R_{s} X \lo R_{s-1} X$. Then each $\Omega_s$ has finite homotopy groups,
 $\omega_s$ is a fibration and there is the induced map  $$K^\land_p=|R_\infty X|\to  \Omega=\varprojlim \Omega_s.$$
Since the homotopy  spectral sequences of both towers converge to 
the homotopy groups of corresponding inverse limits $\pi_*(\varprojlim)$,
 the isomorphism between
  the spectral sequences of the tower of fibrations of simplicial stes $\{R_s X\}$ and its geometric realization $\{\Omega_s\}$
induces an isomorphism of homotopy groups $$\pi_i(K^\land_p)\to  \pi_i(\varprojlim \Omega_s).$$
Hence, the map $K^\land_p\to  \Omega$  is a weak homotopy equivalence. 

The $G$-setting required in the proposition follows from the functoriality of 
 the constructions involved. $\black$
 \\


\

\begin{proposition}
\label{factor}
Let $G$ be a $p$-group,  $K$ a finite $G$-CW-complex and $L$ a  simplicial $G$-complex such that
 all the points  of $L$ have  the same isotropy group.
Given a $G$-map  $f : L\times EG \lo K^\land_p$  where $L \times EG$ is considered 
with the diagonal action of $G$, there is a $G$-map $\phi: L \lo K^\land_p$ such that $\phi\circ\psi$ is $G$-homotopic  to $f$ where 
$\psi:L \times EG\to L$ is the projection.
\end{proposition}
{\bf Proof.} Let $H$ be the isotropy group of the points of $L$. Fix a $0$-simplex (vertex)  $v\in L$
and consider $v\times EG$  and $K^\land_p$  as  $H$-CW-complexes.  
The map $$F_v:v \to {\rm Map}_G (EG, K^\land_p)$$ defined by the restriction $f|_{v\times EG}$ belongs to ${\rm Map}_H(EG, K^\land_p)$.
Indeed, $$F_v(hy)=f(v,hy)=f(hv,hy)=hf(v,y)$$ for all $h\in H$ and $y\in EG$.
By Theorem~\ref{EqSuCo} the inclusion of constant maps 
$$(K^\land_p)^H\stackrel{\subset}\to {\rm Map}_H(EG, K^\land_p)$$  is a weak homotopy equivalence. 
Thus, the map $$F_v:v  \to {\rm Map}_H(EG, K^\land_p)$$ can be deformed to $(K^\land_p)^H$. 
Such deformation  defines
an $H$-homotopy  $$T_v: v\times EG \times [0,1]\lo K^\land_p$$ such that
$T_v (v,  y,0)=f(v,y)$ and $T_v$ is constant on $v \times EG\times 1$. We extend $T_v$  
over $Gv \times EG$
to the homotopy $T_{Gv}$ by $$T_{Gv}(gv, gy, t)=gT_v(v, y, t)\ \ \text{for}\ \ {g\in G},\  y\in EG.$$
 Note that $T_{Gv}$ is
a well-defined $G$-homotopy on $Gv  \times EG$,  since
  for $h \in H$ we have  $hv=v$ and hence
$$T_v (hv,hy, t)=T_v(v,hy,t) =hT_v(v,y,t).$$ 
By performing this procedure independently for all the orbits  of
$0$-simplexes of $L$ we obtain a $G$-homotopy $T_{L^{(0)}}$ that $G$-homotopes 
$f$ restricted to $L^{(0)} \times EG$  to a $G$-map  which  is  constant on $v\times EG$ for all $v\in L^{(0)}$.
Thus we  can replace  $f$ by a $G$-homotopic map  with the  restriction to $v\times EG$ being constant for each $v\in L^{(0)}$.

 Now assume that the restriction of $f$ to $x\times EG$ is constant for all $x\in L^{(n)}$ in the $n$-skeleton  of $L$.
Fix an $(n+1)$-simplex $\Delta$ of $L$. Similarly to what we did before, consider 
$\Delta\times EG$ and $K^\land_p$ as  $H$-CW-complexes. Let $$F_\Delta:\Delta \to {\rm Map}_H(EG, K^\land_p)$$ 
be the map defined by the restriction of $f$ to $\Delta\times EG$.  Note that $F_{\Delta}(\partial \Delta)\subset (K^\land_p)^H$.
Since, in view of Theorem~\ref{EqSuCo},$$\pi_{n+1}( {\rm Map}_H(EG, K^\land_p),(K^\land_p)^H)=0,$$
 there is a ${\rm rel }\ \partial \Delta$ deformation of $F_\Delta$ to $(K^\land_p)^H$. This deformation defines
 an $H$-homotopy $$T_\Delta:\Delta\times EG\times [0,1] \to K^\land_p$$
relative to  $(\partial \Delta) \times EG$ such that
$T_\Delta (x,y,0)=f(x,y)$ and $T_\Delta(x,y,1)$ is constant on $x \times EG$ for all $x\in\Delta$.
We extend $T_\Delta$ to a $G$-homotopy  $T_{G\Delta}$ on $G\Delta \times EG$   by defining
$$T_{G\Delta}(gx, gy, t)=gT_\Delta(x,y, t),\ \ g\in G,\ \ x \in \Delta,\ \ y \in EG.$$ We perform
this procedure independently for all the orbits of $(n+1)$-simplexes of $L$ and this way  we
$G$-homotope $f$ on  $ L^{(n+1)}\times EG$ to a $G$-map which is 
constant on $x\times EG$ for all $x\in L^{(n+1)}$. Thus we can $G$-homotope $f$ to a $G$-map which is
constant on $x\times EG$ for all $x\in L$  and the proposition follows. $\black$
\\
\\
In the proof of Proposition \ref{factor} we showed that the projection of $L \times EG$ to $L$ 
induces a surjection from  the path connected components of  ${\rm Map}_G(L, K^\land_p)$ to
the path  connected components of ${\rm Map}_G(L \times EG, K^\land_p)$.   We note that
 Proposition \ref{factor} can be generalized as follows:
The restriction on the isotropy groups of the points of $L$ can be dropped and 
it can be shown that the projection of $L \times EG$ to $L$ induces a weak homotopy equivalence 
 $${\rm Map}_G(L, K^\land_p) \to {\rm Map}_G(L \times EG, K^\land_p).$$

\end{section}

\section{G-maps and G-homotopies}
\label{g-maps}

In this section we extend two well-known facts in topology to the category $G$-spaces and $G$-maps for  a discrete group $G$.
The first fact is the following
\begin{proposition}\label{lifting}
Let $\phi: X\lo Y$ be a  weak homotopy equivalence  between path connected spaces and let $L$ be a  simplicial complex. Then
every map $f:L\to Y$ admits a homotopy lift $\bar f:L\to X$.
\end{proposition}
The proof goes like this. Replacing $Y$ by  the mapping cylinder of $\phi$,
 we may assume without loss of generality that $\phi$ is an embedding. 
 The weak homotopy equivalence condition implies that the relative homotopy groups $\pi_i(Y,X)=0$. 
 By induction on $n$ we can deform the restriction of $f$ to the $n$-skeleton $L^{(n)}$ to $X$
 by a ${\rm rel } \ L^{(n-1)}$ deformation  for $n>0$. This defines a homotopy lift of $f$.

We recall that for the maps $\phi:X\to Y$ and $f:L\to Y$, a map $\bar f:L\to X$ is a homotopy lift of $f$ with respect to $\phi$ if $\phi\circ \bar f$ is homotopic to $f$.
In the case of $G$-spaces and $G$-maps we define a $G$-homotopy lift  $\bar f$ by imposing the requirement that $\phi \circ \bar f$ is $G$-homotopic to $f$.

Then the $G$-extension of the first fact is
\begin{proposition}
\label{weak-equivalence-lifting}
Let $\phi: X\lo Y$ be a $G$-map of  path connected topological $G$-spaces such that
$\phi$ is a weak homotopy equivalence and let $L$ be a free simplicial $G$-complex. Then
every $G$-map $f:L\to Y$ admits a $G$-homotopy lift $\bar f:L\to X$.
\end{proposition}
{\bf Proof.} Again by taking  the mapping cylinder we may assume that $\phi$ is an embedding. Thus, $\pi_i(Y,X)=0$.
As above by induction on $n$ we can deform the restriction  $f|_{L^{(n)}}$ to a map with the image in $X$
by a ${\rm rel} \ L^{(n-1)}$ deformation.
Since the $G$-action on $L$ is free, we can do it in the equivariant setting. For that for each $n$ we consider separately the orbits of all $n$-simplices.
Do such deformation for an arbitrary representative of the orbit and extend it to the other elements by the $G$-translation.
$\black$
\\

The second fact is slightly more delicate. 
Let $$\Omega_1 \stackrel{\omega^2_1}\leftarrow \Omega_2\stackrel{\omega^3_2}\leftarrow \Omega_3\leftarrow\cdots $$ be a tower of fibrations with 
$\Omega=\varprojlim \Omega_i$ and let $L$ be a CW-complex which is the union of an increasing sequence of subcomplexes $L_i, i\geq1$.
We define the function $$\psi_i^{i+1}:[L_{i+1},\Omega_{i+1}]\to[L_{i},\Omega_{i}]$$ 
as $$\psi^{i+1}_i([f])=[(\omega^{i+1}_i \circ f)|_{L_i}]$$ where $[X,Y]$ denotes the homotopy classes of maps $f:X\to Y$.
Then the projection  $\omega^\infty_i:\Omega\to\Omega_i$ of the restriction to $L_i$ defines 
the  function of homotopy classes $$\Phi_i:[L,\Omega]\to[L_i,\Omega_i].$$

\begin{proposition}\label{surjective}
The function $\Phi=\varprojlim \Phi_i:[L,\Omega]\to \varprojlim \{[L_i,\Omega_i],\psi^{i+1}_i\}$ is surjective.
\end{proposition}
Proposition~\ref{surjective} says that given a sequence of maps 
$f_i:L_i\to\Omega_i$ such that $(\omega^{i+1}_i\circ f_{i+1})|_{L_i}$ is homotopic to $f_i$ we can find a map $f:L\to\Omega$
such that $f_i$ is homotopic to $\omega^\infty_i \circ f|_{L_i}$ for all $i$. Such map can be constructed by
a recursive application of the Homotopy Lifting Property for pairs.
We recall that a map $\phi : E\to B$  satisfies the {\em Homotopy Lifting Property for a pair} $(X,A)$ if for any
homotopy $H:X\times [0,1]\to B$ and any lift $\bar H' : X' \lo E$ of $H$ restricted to 
$X'=X \times 0 \cup A \times [0,1]$ 
there is a lift $\bar H:X\times [0,1]\to E$
of $H$ which extends  $\bar H'$.
The following is well-known~\cite{H}:
\begin{theorem}\label{lift}
Any Hurewicz fibration $\phi :E\to B$  satisfies the Homotopy Lifting Property for CW-complex pairs $(X,A)$.
\end{theorem}

The $G$-version of Proposition~\ref{surjective} is the following
\begin{proposition}\label{G-version surjective}
Suppose $\{ \Omega_i \}$ is a tower of fibrations of $G$-spaces whose bonding maps are $G$-maps and
let $L$ be a free simplicial $G$-complex which is the union of an increasing sequence of 
$G$-subcomplexes $L_i$.
Then the function
$$\Phi_G:[L,\Omega]_G\to\varprojlim [L_i,\Omega_i]_G$$ is surjective.
\end{proposition}
Here  for $G$-spaces we denote by $[X,Y]_G$  the set of $G$-homotopy classes of $G$-maps from $X$ to $Y$.
The functions     $\Phi_i$ and $\Phi$  defined before can be translated  in an obvious way to the $G$-setting to the functions
  $$\Phi_i^G:[L,\Omega]_G\to[L_i,\Omega_i]_G\ \ \text{and}\ \  
 \Phi_G=\varprojlim \Phi_i^G.$$
 We note that   Propositions \ref{surjective} and \ref{G-version surjective} can be strengthened  by showing that the functions
 $\Phi$ and $\Phi_G$ are bijections provided  all $L_i$ are finite sucomplexes and $\Omega_i$ have finite homotopy groups.

The proof of Proposition~\ref{G-version surjective} is a recursive application of the following $G$-version of the Homotopy Lifting Property for pairs.
\begin{theorem}
\label{fibration}
Let $E$ and $B$ be $G$-spaces,  $\phi : E \lo B$  a fibration which is a $G$-map and
let $L$ be a free simplicial $G$-complex, $A$ a $G$-subcomplex of $L$,
and $H: L \times [0,1] \lo B$ a  $G$-homotopy. Let $L'=L \times 0 \cup A \times [0,1]$ and
 let $\bar H':  L' \to E$ be a $G$-lift of  $H$ restricted to ${L'}$. 
Then there is a $G$-lift $\bar H:L\times [0,1]\to E$ of $H$ that extends $\bar H'$.
\end{theorem}
{\bf Proof.}  We define $$M_n =(L\times 0)\cup (A \times [0,1]) \cup (L^{(n)}\times [0,1]).$$ 
First we show that  $\bar H'$ extends to a $G$-map $\bar H_0 : M_0 \lo E$ which is a lift of  $H$ restricted to $M_0$.
Take a $0$-simplex $v\in L\setminus A.$   Since $\phi$ is a fibration we can lift $H|_{v \times [0,1]}$
to a path $$\bar H_v : v \times [0,1]\to E$$  such that $\bar H_v(v,0)=\bar H'(v)$.
 We extend $\bar H_v$ over $Gv \times [0,1]$ to
a $G$-homotopy $\bar H_{Gv} : Gv \times [0,1] \lo E$ by 
$$\bar H_{Gv} (ga, t)=g{\bar H_{v}}(a,t), \ \ {g\in G,}\ \ a \in \Delta,\ \ t\in [0,1].$$ We do this independently for
all the orbits of $0$-simplexes of $L\setminus A$ to extend $\bar H'$ to
a $G$-map $\bar H_0$.

Now assume that we already {have} constructed 
 a $G$-map $\bar H_n : M_n \lo E$ which is a lift of  $H$ restricted to $M_n$ and extends $\bar H_{n-1}$.
 Consider an $(n+1)$-simplex $\Delta$ of $L$ not contained in $A$. 
 Since $\phi$ is a fibration, we can apply the homotopy lifting property 
  for the pair $(\Delta,  \partial \Delta )$
 to lift $H$ restricted to $\Delta \times [0,1]$ to
 a map $$\bar H_\Delta : \Delta \times [0,1] \lo E$$ that extends $\bar H_n$ restricted to $ (\Delta\times [0,1])\cap M_{n}$.
 Now extend $\bar H_\Delta$ to a $G$-homotopy
 $$\bar H_{G\Delta}: G\Delta \times [0,1] \lo E$$ by the formula $$\bar H_{G\Delta} (ga, t)=g{\bar H}_\Delta (a,t),\ \ {g\in G,}\ \ a \in \Delta.$$
 Doing that independently for  all the orbits of $(n+1)$-simplexes not contained in $A$ we get
 a $G$-map $$\bar H_{n+1} : M_{n+1} \lo E$$ which extends $\bar H_n$ and which is a lift
 of $H$ restricted to $M_{n+1}$.
 
 Then the $G$-maps $\bar H^n$ for all $n$ provide the $G$-lift $\bar H$ required in the proposition.
 $\black$

\begin{proposition}\label{finite G-sets}
Let  $L$ be a finite simplicial free $G$-complex. Then for any path connected $G$-space  $X$ with finite homotopy groups the set $[L, X]_G$  is a finite.
\end{proposition}
{\bf Proof}.
We use  induction on $n=\dim L$. The statement is trivial for $n=0$. Suppose it holds true for $n-1$ and $\dim L=n$. 
Let $f:L^{(n-1)}\to X$ and let $\Delta$ be an $n$-simplex in $Y$.
Since $\pi_n(X)$ is finite there are finitely many non-homotopic maps $f_\Delta:\Delta\to X$ extending $f|_{\partial\Delta}$. 
Then there are finitely many $G$-homotopy classes
$[L^{(n-1)}\cup G\Delta, X]_G$. We apply induction on the number of the $G$-orbits of $n$-simplices in $L$ to complete the proof.
$\black$

\begin{proposition}
\label{cantor}
Let $\{ \Omega_i\}$ be a tower of fibrations of path connected    $G$-space with finite homotopy groups 
and bonding maps
$\omega_i : \Omega_i \lo \Omega_{i-1}$ being $G$-maps  and fibrations  and let 
$\Omega=\varprojlim \Omega_i$. Suppose that $L$ is a free simplicial $G$-space such that
$L$ is the union of an increasing sequence of finite  $G$-subcomplexes $L_i$  and
each $L_i$ admits a $G$-map to $\Omega$. Then $L$ also admits a $G$-map to $\Omega$.
\end{proposition}
{\bf Proof}. Since each $L_i$ admits a $G$-map to $\Omega$, we obtain $[L_i,\Omega_i]_G\ne\emptyset$. By Proposition~\ref{finite G-sets} the sets 
$[L_i,\Omega_i]_G$ are finite. Therefore, $$\varprojlim[L_i, \Omega_i]_G\ne\emptyset.$$
Proposition~\ref{G-version surjective} implies that $[L,\Omega]_G\ne \emptyset$.
$\black$\\

We will summarize the results of this section in

\begin{proposition}\label{summary}
Let $K$ be a connected finite $G$-CW-complex and $L$  a free simplicial $G$-complex such that 
$L$ is  the union of an increasing sequence of finite  $G$-subcomplexes $L_i$ 
such that each $L_i$ admits a $G$-map to $K^\land_p$. Then $L$ admits a $G$-map to $K^\land_p$.
\end{proposition}
{\bf Proof.}  Consider the tower of fibrations $\{ \Omega_s \}$ satisfying the conclusions of
Proposition \ref{tower}  and let  $\phi : K^\land_p \lo \Omega=\varprojlim \Omega_s$
be a $G$-map which is also a weak homotopy equivalence. Then each $L_i$ admits a $G$-map
to $\Omega$ and, by Proposition  \ref{cantor}, we get that $L$ admits a $G$-map to $\Omega$.
By Proposition \ref{weak-equivalence-lifting}, the latter map admits a $G$-homotopy lift
to $K^\land_p$. $\black$

\begin{section}{Proof of Proposition \ref{proposition-for-borsuk-ulam} and  generalizations of 
Theorems A and B}
\label{section-borsuk-ulam}

{\bf Proof of Proposition \ref{proposition-for-borsuk-ulam}.} 
We recall  that    $EG$ is a  contractible  $G$-CW-complex with a free action of $G$.
Thus if $G$ is  a finite cyclic group
we can represent $EG$ as an increasing sequence of spheres $\s^{2n-1}$ 
equipped with standard action of $G$  such that $\s^{2n-1}$ is a $G$-subcomplex 
of $\s^{2n+1}$.  We will refer to this model of $EG$ as the infinite dimensional sphere $\s^\infty$.

By Remark  \ref{remark},  we can assume that
all the spaces
$K(l_q,m, q)$ considered in this proposition are connected.

We prove the proposition by induction on $m$ and start with $m=1$.
Set $N_0$ to  be a singleton with 
the  trivial action of $G_1$. {We} recall that    $N_1=\s^{2n -1}$  for $n=n_1$ and represent $N_1$ as
 $N_1=N_0 \times \s^{2n -1}$ with the diagonal action
of $G_1$. We need to find $n_1$ satisfying the conclusion of the proposition for $N_1$.
Assume that such $n_1$  does not exist. This means 
 that  $N_0 \times \s^{2n -1}$  does admit 
a $G_1$-map to  $K(l_0,1,0)^\land_p$ for every  $n$.

Represent  $EG_1$ as the infinite dimensional sphere $\s^\infty$.
Then, by Propositions \ref{summary}, 
 $N_0 \times EG_1$  admits a $G_1$-map  $$f:N_0\times EG_1\to K(l_0,1,0)^\land_p$$ and hence, by 
Proposition \ref{factor},  $N_0$ admits a $G_1$-map  $$\phi:N_0\to K(l_0,1,0)^\land_p$$  as well. Recall that
$K(l_0,1,0)$ does not have  fixed points of the $G_1$-action. Then, by Proposition \ref{basic-1},
$K(l_0,1,0)^\land_p$ does not have fixed points either and this contradicts the fact
that $N_0$  admits a $G_1$-map to $K(l_0,1,0)^\land_p$.

Assume that we already determined $n_1,\dots,n_m$ satisfying the conclusion of the proposition
and proceed to $n_{m+1}$. Again assume that $n_{m+1}$  with required properties
does not exist and  recall that $N_{m+1}=N_m \times \s^{2n-1}$ with $n=n_{m+1}$
and the diagonal action of $G_{m+1}$.
Thus we assume that  there is $q\leq  m$ such that  $N_m \times \s^{2n-1}$ admits a $G_{m+1}$-map 
$$f_n:N_m \times \s^{2n-1}\to K(l_q, m+1, q)^\land_p$$ for every $n$.
Then,
representing  $EG_{m+1}$ as the infinite dimensional sphere
 $\s^\infty$  and
 applying  Propositions \ref{summary}, we conclude that $N_m \times EG_{m+1}$
 admits a $G_{m+1}$-map $$f:N_m \times EG_{m+1}\to
 K(l_q, m+1, q)^\land_p.$$
   Consider    $N_m$   as a $G_{m+1}$-complex with the action of $G_m$ extended 
to   the action of $G_{m+1}$ whose isotropy group is $H=\z/p\z$.
Then, by Proposition \ref{factor}, we get that $N_m$
admits a $G_{m+1}$-map  $$ \phi:N_m\to K(l_q,m+1,q)^\land_p.$$  Since all the points
of $N_m$ are $p^m$-periodic,  the map $\phi$ lands in $(K(l_q, m+1, q)^\land_p)^H$ and
hence we have a $G_{m+1}$-map  $$\phi:N_m\to (K(l_q, m+1, q)^\land_p)^H.$$

Let $q=m$. {We} recall that $K(l_{m},m+1,m)$ does not have $p^m$-periodic points.
 Then Proposition \ref{basic-1}
implies that $K(l_{m},m+1,m)^\land_p$ does not have $p^m$-periodic points either.
This means that $$(K(l_m, m+1, m)^\land_p)^H=\emptyset$$  
 and we get a  contradiction.
 
Now  let $q < m$.  By Proposition \ref{basic-1}, $N_m$ admits a $G_{m+1}$-map 
 $$\psi :N_m\to(K(l_q, m+1, q)^H)^\land_p.$$ By Proposition \ref{general-3} 
 we have a $G_m$-isomorphism $$K(l_q, m+1, q)^H\cong K^m(l_q, m, q)$$  
 where the first complex is considered with the action of $G_{m+1}/H=G_m$. Then the functoriality
of the $p$-adic completion  implies a $G_m$-isomorphism  $$(K(l_q, m+1, q)^H)^\land_p\cong K(l_q, m, q)^\land_p.$$  Since $H$ acts trivially on $N_m$, 
	and $G_m=G_{m+1}/H$ we get that $\psi$ defines
	 a $G_m$-map  $$N_m\to K(l_q, m, q)^\land_p$$ which brings a contradiction.
$\black$
\\
\\
Theorem A can be generalized as follows.

\begin{theorem}
\label{generalization-1}
Let $p$ be a prime number and let $A$ be an infinite  subgroup of a countable  product  of finite $p$-groups.
We consider $A$ as a discrete  group. Then $A$ admits a free action by isometries on a compact metric space $X$
such that the dynamical system $(X,A)$ is not equivariantly embeddable into the cubical shift 
$([0,1]^k)^A$  for every $k$.
\end{theorem}
Note that Theorem \ref{generalization-1} covers a relatively large variety of groups $A$,  e.g.
 $\z^n$,  $\z_p$, $(\z/p^n\z)^\N$.
 Theorem \ref{generalization-1} is derived from the following  Borsuk-Ulam  type theorem.
\begin{theorem}
\label{generalization-2}
 Let $p$ be a prime number and let  $G$ be the product of non-trivial  finite $p$-groups 
$H_i, i\geq 1$. Set $H_0$ to be the trivial group. 
Consider  each group $H_i$ as a discrete group  and the group $G$ as 
a $0$-dimensional compact metric group with 
the product topology. Then for every sequence of natural numbers $l_q, q\geq 0,$ there is 
a free action of $G$ on a compact metric space $X$ such for every $q\geq 0$  and every
map $f : X \lo \re^{l_q}$ there is a point $x \in X$ such that 
$f^{G}(Gx)$ contains  at most  $|H_0 \times \dots \times H_q|$   points 
where
$f^G : X \lo (\re^{l_q})^G$ is the  $G$-map induced by $f$. 
In particular, $f(Gx)$ contains at most $|H_0 \times \dots\times H_q|$  points.
\end{theorem}
{\bf Proof of Theorem \ref{generalization-1}.}  Let $X$ be as in Theorem \ref{generalization-2}.
Since  $G$ is compact we may assume without loss of generality that
$G$ acts on $X$ by isometries \cite{montgomery-zippin}. 
Note that for a subgroup $A$ of $G$, $x \in X$ and a map
 $f: X \lo \re^{l_q}$
we have that $$|f^A(Ax)|\leq |f^G(Gx)|$$ where the $A$-map  and 
the $G$-map $$f^A : X \lo (\re^{l_q)})^A\ \ \text{and}\ \ f^G : X \lo (\re^{l_q)})^G $$ are  
induced by $f$.
Then the proof of Theorem A  applies for  the sequence $l_q=q+1$
to  get the required result. $\black$\\

  {Below we} outline how the approach of this paper can be adjusted  to prove Theorem \ref{generalization-2}.
  For a finite group  $H$ represent $EH$ (a contractible $H$-CW-complex with a free action of $H$) as 
  the  union of an increasing sequence of  finite  simplicial  $H$-subcomplexes  $EH(n), n\geq 1$.
  For a cyclic group $H$ we can consider the infinite sphere $\s^\infty$ as a model for $EH$ and
  set $EH(n)=\s^{2n-1}$ with the standard action of $H$. For a general finite group $H$ we can consider
 Milnor's model  of $EH$ with $EH(n)=H*\dots * H$  ($n$ factors) and the  action of $H$ 
 on $EH(n)$ induced by the action of $H$ on itself.
 
  Let $H_i, i\geq 1,$ be a sequence of non-trivial finite $p$-groups, set  
  $H_0$ to be the trivial group and  denote 
  $$G_j=H_0 \times \dots \times H_j,\ j\geq 0.$$Clearly 
    $H_j$,  $G_j, j\leq m, $  and $$G_m/G_j=H_{j+1} \times \dots \times H_m,\ j<m$$   can be considered 
      as subgroups of $G_m$.
  We adjust the notations of Section \ref{reduction}
  as follows: $E(l,m)=(\re^l)^{G_m}$ with the shift action of $G_m$, 
  $$E^q(l,m)=(E(l,m))^{G_m/G_q}\ \ \text{and}\ \ E(l,m,q)=E(l,m)\setminus E^q(l,m)\ \text{for}\ 0\leq q< m.$$
   The   $G_m$-complex $K(l, m, q)$ is defined  as in Section \ref{reduction}. 
   The complex  $K^m(l,m+1,q)$   is defined as
  $$K^m(l,m+1,q)=K(l,m+1,q)^{H_{m+1}}$$  and $K^m(l,m+1,q)$ can be considered with 
  the action of  the quotient group $G_m=G_{m+1}/H_{m+1}$.

  Given a sequence of natural numbers $n_i, i\geq 1,$ we define 
  the free simplicial $G_m$-complex $N_m$  as $$N_m=EH_1(n_1) \times \dots \times EH_m(n_m)$$
  with the product action of $G_m$ and define the compact metric $G$-space $X$ as the product
  of  $EH_i(n_i)$ for all $i\geq 1$ with the product action of the compact  group $G$ which  is the product of  all   $H_i, i\geq 1,$
  considered as discrete groups.
  Then all the results of Section \ref{reduction}  hold with obvious adjustments as replacing
  $p^q$ by $|G_q|$,   diagonal actions 
  by  product actions, $\z$ by $G$   etc.
  The proof of Proposition \ref{proposition-for-borsuk-ulam} also works with obvious adjustments.

  Thus we get that that for every sequence $l_q, q\geq 0,$ there is a sequence $n_i, i \geq 1,$
  such that $(X,G)$ constructed as above satisfies  the conclusions of 
   Theorem \ref{generalization-2}. \\

   Let us finally point out  that given a collection of  $p$-groups $H_i$
 our approach  does not
bring an estimate on the growth of  $n_i$,   it is just about  pure existence.  
 In this connection  we mention  a few  related constructive  results obtained by translating 
 Borsuk-Ulam theorems to  a question about the existence of  non-vanishing sections of  certain induced
 vector bundles. 
  Dranishnikov    \cite{dranishnikov-1, dranishnikov-2} showed that 
   if $X$ is   the product of the spheres $\s^{n_i}$  with $n_i=2^i, i\geq 1, $   and $X$ is considered with
   the product action of $G=(\z/2\z)^\N$ then for every map $f : X \lo \re $ there is 
   $x\in X$ such that $f(Gx)$ is a singleton. Dranishnikov's result was extended by Turygin \cite{turygin2}
 for $G=(\z/p \z)^\N$ for all primes $p$ and maps $f:X\to\mathbb R^l$  providing an explicit estimate on the growth of $n_i$.
    Dranishnikov's result  relies on computation of the corresponding Stiefel-Whitney class
    and Turygin's on computation of the mod $p$ Euler class. The Borsuk-Ulam theorem
    for a free action of  $G=\z/p^k\z$ on a sphere $\s^{2n-1}$ and maps to $\re^l$ with an explicit estimate for $n$ depending
    on $p, k$ and $l$ was proved by  Munkholm \cite{munkholm} for odd $p$.
    This  is done by computing 
       the corresponding complex $K$-theory  Euler class.  
       Later Munkholm's result was generalized and extended by the school of Puppe~\cite{Ba} to all primes $p$.    
       Unfortunately    all these methods do not work for the product  action of 
 $$G=\z/p\z\times\z/p^2\z\times\cdots\times \z/p^k\z\times\cdots$$ on the product of spheres. This case is covered
 by Theorem \ref{generalization-2} and demonstrates
    a crucial advantage of our approach based on
   $p$-adic completion and the Sullivan conjecture.

\end{section}

 dranish@math.ufl.edu \\
Department of Mathematics,  University of Florida
444 Little Hall, Gainesville, FL 32611-8105, USA
\\\\
mlevine@math.bgu.ac.il\\
Department of Mathematics, Ben-Gurion University of the Negev
P.O.B. 653, Be’er Sheva 84105, Israel

\end{document}